\documentclass{amsart}
\usepackage{amsfonts}
\usepackage{amssymb}
\usepackage{amsmath}
\usepackage{amsthm}
\theoremstyle{plain}
\newtheorem{thm}{Theorem}[section]
\newtheorem{lemma}[thm]{Lemma}
\newtheorem{prop}[thm]{Proposition}

\theoremstyle{definition}
\newtheorem{definition}[thm]{Definition}
\newtheorem*{ex}{Example}
\numberwithin{equation}{section}
\begin{document}
\title{Parameter Families of Iterated Function Systems and Continuity}
\author{Max Murphy}
\address{Department of Mathematics and Statistics\\University of
Nevada, Reno}
\thanks{This project was funded by the Office of Undergraduate Research, University of Nevada, Reno; and by the NSF under grant number 0447416. Further thanks go to my faculty mentor, Dr. Alex Kumjian.}
\date{\today}
\begin{abstract}In this paper, we investigate parameter families of iterated function systems and continuity. Specifically, if we have a set of iterated function systems that depend continuously on a parameter, which properties of the invariant sets will vary continuously? We show here that the invariant sets will vary continuously as a function of the parameter. Furthermore, we give some sufficient conditions under which the Hausdorff dimension will depend continuously on the parameter. Lastly, we give a counter-example where the dimension is not a continuous function of the parameter.
\end{abstract}
\maketitle
\section*{Introduction}
Here we study parameter families of iterated function systems and
continuity. If we have a parameter family of iterated function
systems, then we can ask questions about the sets generated by the
iterated function system, and how they relate to the parameter. This
is important to physical applications as well as abstract mathematics,
because if we make a ``good'' measurement of the parameter, we would
hope that we will be able to draw ``good'' conclusions.  For
example, since every iterated function system generates an invariant
set, we can ask whether these invariant sets depend continuously on the parameter. Furthermore, since any set has a
dimension associated with it, we can ask whether the dimension of
the invariant sets vary continuously. \\
\indent First, we go over some preliminary definitions and theorems
concerning fractal geometry, particularly iterated function systems.
For a more detailed introduction to fractal geometry and iterated
function systems, see \cite{MTFG} or \cite{FG}. After reviewing some
of the basics, we then go into the results that concern continuity
of parameter
families of iterated function systems. A parameter family of IFS's is basically defined to be a family of iterated function systems that depend continuously on some parameter. Thus, instead of generating a single fractal with a single IFS, we generate an entire family of fractals; one for each parameter value. We then explore how the properties of these fractals change with the parameter, and what properties they share. \\
\indent For our first result, we show that the invariant sets
(equipped with the Hausdorff metric) do vary continuously as a
function of the parameter. After showing that the invariant sets
themselves work nicely with respect to continuity, we move on to the
Hausdorff dimension. Unfortunately, the Hausdorff dimension will not
necessarily be continuous as a function of the parameter. However,
if we impose the right hypotheses on the iterated function system,
we can show that the dimension will be continuous. Specifically, if
the family of iterated function systems consist of similarities that
satisfy the open set condition, we prove that the dimension must be
continuous. We then give a counterexample to show that the dimension
can have discontinuities if the hypotheses are not satisfied.
\section{Preliminaries}
Before we can showcase our results, we first give some basic definitions and theorems concerning iterated function systems, which we use in the proofs below. Most of the notation here is taken from \cite{MTFG}.
\begin{definition}
Given two metric spaces $X$ and $Y$, we call a function $f:X \to Y$ a \emph{similarity} if, for some $r \in (0,\infty)$, \[d(a,b) = r\,d(f(a),f(b))\] for all $a,b \in X$.
\end{definition}
\begin{definition}
On a metric space $X$, a \emph{contraction map} is a function $f:X
\to X$ such that, for some $r \in (0,1)$, \[d(f(x),f(y)) \leqslant r
\, d(x,y)\] for all $x,y \in X$. Note that every contraction map is
continuous.
\end{definition}
One of the most important facts about contraction maps, and one that we will repeatedly use, is the fixed point theorem for contraction maps. It shows that any contraction map will have a point which is fixed, and furthermore, it gives us a way to actually find this point.
\begin{thm}\cite[p. 50]{MTFG}
Let $X$ be a complete metric space, and let $f:X \to X$ be a
contraction map. Then $f$ has a unique fixed point, that is, there
is a unique $x \in X$ such that $f(x) = x$. Furthermore, if
$\{x_n\}$ is a sequence that satisfies $x_{n+1} = f(x_n)$ (where
$x_0$ is an arbitrary point in $X$), then $x = \lim x_n$.
\end{thm}
Now, we need to define the Hausdorff metric, which is a way to measure the distance between the compact sets in any metric space. This is a very powerful tool, because it will give us a way to compare the fractals we have constructed.
\begin{definition}
Given a metric space $X$ and a subset $A$, we define $N_r(A)$, \emph{the open r-neighborhood of A} to be \[N_r(A) = \{x \in X : \exists \:a \in A \; \text{such that} \;\; d(a,x) < r\}.\]
\end{definition}
\begin{definition}
Given a metric space $X$, we define $\mathfrak{K}(X)$ to be the set
of all compact subsets of $X$. This set is usually equipped with the
\emph{Hausdorff metric}, $h$. where \[ h(A,B) = \inf\{r : A
\subseteq N_r(B) \; \text{and} \; B \subseteq N_r(A)\}\] It isn't too hard to show that this is indeed a metric \cite[p. 66]{MTFG}.
\end{definition}
Equivalently, the Hausdorff distance between $A$ and $B$ is the smallest number $r$ such that every point in $A$ is within $r$ of some point in $B$, and every point in $B$ is within $r$ of some point in $A$.
\begin{prop}\cite[p. 67]{MTFG}
Let $X$ be a complete metric space. Then the space $\mathfrak{K}(X)$ is also complete.
\end{prop}
The fact that $\mathfrak{K}(X)$ is also complete is very important, because we can now apply the fixed point theorem for contraction maps to that space. These sets that are fixed points are often fractals.
\begin{definition}
Given a complete metric space $X$, we define an \emph{iterated
function system} to be a finite number of contraction maps $f_i:X
\to X$ for $i \in \{1,2,...n\}$.
\end{definition}
\begin{definition}
Let $X$ be a complete metric space and suppose we are given an
iterated function system, $\{f_i\}$ on $X$. We define the
\emph{invariant set} of the iterated function system to be a fixed
point of $F: \mathfrak{K}(X) \to \mathfrak{K}(X)$, where \[F(C) =
\bigcup_{i=1}^{n}f_i(C).\] Note that this function $F$ is a
contraction map. Since $\mathfrak{K}(X)$ is complete by proposition
1.6, the invariant set of an iterated function system is always
uniquely defined, by theorem 1.3.
\end{definition}
We now move on to the subject of dimension. There are some measure
theoretic ideas going on behind this, but we omit discussing this.
For more detail on measures and fractals, see \cite{IPFM}.
\begin{definition}
Let $X$ be a metric space, and $C \subseteq X$. For any countable family of sets, $\{U_i\}$, we say $\{U_i\}$ is a \emph{$\epsilon$-cover} of $C$ if $C \subseteq \bigcup_{i=1}^{\infty} U_i$ and, for each $i$, $|U_i| < \epsilon,$ where $|U_i|$ denotes the diameter of $U_i$. For any $s \geqslant 0, \epsilon > 0$, we define \[
\mathcal{H}_{\epsilon}^s(C) = \inf \left\{\sum_{i=1}^{\infty}|U_i|^s : \{U_i\} \; \text{is an $\epsilon$-cover of $C$} \right\}.\]
Furthermore, we define $\mathcal{H}^s(C)$, the \emph{s-dimensional Hausdorff outer measure}, to be \[
\mathcal{H}^s(C) = \lim_{\epsilon \to 0} \mathcal{H}_{\epsilon}^s(C) \]
\end{definition}
Intuitively, we think of $\mathcal{H}^s(C)$ to be the s-dimensional
``volume,'' for example, think of $\mathcal{H}^2(C)$ as the area
contained in $C$. However, we have the advantage that $s$ here does
not need to be an integer, so this generalizes to non-integer
dimensions, which is our goal.
\begin{prop}\cite[p. 31]{FG}
Let $X$ be a metric space, and let $C \subseteq X$. There is a
unique value $s_0 \in [0,\infty]$ such that $\mathcal{H}^s(C) = 0$
whenever $s > s_0$, and $\mathcal{H}^s(C) = \infty$ whenever $s <
s_0$. We call this value $s_0$ the \emph{Hausdorff dimension} of
$C$.
\end{prop}
Informally, if a set has Hausdorff dimension $s$, then it has no $t$-dimensional ``volume'' for $t>s$, and it has infinite $t$-dimensional ``volume'' when $t<s$. \\
\indent Of course, Hausdorff dimension is very difficult to
calculate in general. However, there is a theorem which allows us to
calculate the dimension of set in special cases. First though, we
have to make another definition.
\begin{definition}
Given an iterated function system $\{f_i\}$ on a complete metric space $X$, we say that the iterated function system satisfies the \emph{open set condition} if there is some nonempty open set $U$ such that
\[\begin{aligned} &f_i(U) \subseteq U \quad \text{for all $i$} \\
&f_i(U) \cap f_j(U) = \varnothing \quad \text{for all $i \ne j$}
\end{aligned} \]
\end{definition}
\begin{thm}\cite[p. 130]{FG}
Let $\{f_i\}$ be an iterated function system on a complete metric space $X$, consisting entirely of similarities where $f_i$ has similarity ratio $r_i$. Furthermore, suppose that $\{f_i\}$ satisfies the open set condition. Let $A$ be the invariant set of $\{f_i\}$, and let $s$ be the Hausdorff dimension of $A$. Then $s$ is the unique number that satisfies the equation \[
\sum_{i=1}^{n}r_i^s = 1. \]\end{thm}
\section{Results}

For our first theorem, we show that if we have a continuously
parameterized family of iterated function systems, then the
invariant set varies continuously as a function of the parameter. To
do this, we first need to show that the fixed point of a family of
contraction maps depends continuously on the parameter.
\begin{prop}
Let $P$ and $X$ be metric spaces, the latter being complete. Let $f: P \times X \to X$ be a family of contraction maps. That is, $x \mapsto f(p,x)$ is a contraction map with ratio $r_p$ for any $p \in P$. Furthermore, suppose that $p \mapsto f(p,x)$ is a continuous function for any fixed $x \in X$. Then the fixed point $x_p$ depends continuously on $p$.
\end{prop}
\noindent \emph{Proof.} Let $p \in P$ and $\epsilon > 0$ be given. Then, for any $q \in P$
\[\begin{aligned}
d(x_p,x_q) &\leqslant d(x_p,f(p,x_q)) + d(f(p,x_q),x_q) \\
&= d(f(p,x_p),f(p,x_q)) + d(f(p,x_q),f(q,x_q)) \\ &\leqslant r_pd(x_p,x_q) + d(f(p,x_q),f(q,x_q)).
\end{aligned}\]
These relations all follow from the definition of contraction, and the triangle inequality. Thus,
\[
(1-r_p)d(x_p,x_q) \leqslant d(f(p,x_q),f(q,x_q)).
\]
$f(\cdot,x)$ is continuous for any fixed $x \in X$, so there is a $\delta$ such that
\[
d(f(p,x_q),f(q,x_q)) < (1-r_p)\epsilon
\]
whenever $d(p,q) < \delta$. Thus, if $d(p,q) < \delta$, then
\[
d(x_p,x_q) \leqslant \frac{d(f(p,x_q),f(q,x_q))}{1-r_p} < \frac{(1-r_p)\epsilon}{1-r_p} = \epsilon.
\]
Thus, $p \mapsto x_p$ is continuous. $\hfill \square$
\\ \\
\indent Next, we just need to show that the entire iterated function system depends continuously on the parameter in $\mathfrak{K}(X)$, and then the result will follow.
\begin{lemma}
Let $P$ and $X$ be metric spaces, the former being compact, the latter being complete. Let $f_i: P \times X \to X$ be continuous for each $i \in \{1,2,...,n\}$ Let $F: P \times \mathfrak{K}(X) \to \mathfrak{K}(X)$ be defined as $F(p,C) = \bigcup_{i=1}^{n}f_i(p,C)$. Then $F(\cdot,C)$ is continuous for any $C \in \mathfrak{K}(X)$
\end{lemma}
\noindent \emph{Proof.} Let $\epsilon > 0.$ Let $C$ be a compact
subset of $X$. $f_i$ is continuous on $P \times C$, and since this
is a compact space, it is uniformly continuous there. Thus there is
a $\delta$ such that $d(f_i(p,x),f_i(q,y)) < \epsilon$ whenever
$\max(d(p,q), d(x,y)) < \delta$. Suppose $d(p,q) < \delta$. Then,
for all $c \in C, \: d(f_i(p,c),f_i(q,c)) < \epsilon$. Thus, by the
definition of the Hausdorff metric, $h(f_i(p,C),f_i(q,C)) <
\epsilon$. And again from the definition $h(F(p,C),F(q,C)) <
\epsilon$. Thus $p \mapsto F(p,C)$ is continuous for any $C \in
\mathfrak{K}(X)$. $\hfill \square$

\begin{thm}
Let $P$ and $X$ be metric spaces, where $P$ is compact and $X$ is
complete. Let $f_i:P \times X \to X$ be a continous family of
iterated function systems. Let $A(p)$ be the invariant set of the
system for any $p \in P$. Then $p \mapsto A(p)$ is a continuous
function.
\end{thm}
\noindent \emph{Proof.} By lemma 2.2, the function $F$ depends continuously on $p$ for any fixed $C$. By proposition 2.1, its fixed point depends continuously on $p$. $\hfill \square$
\\ \\
\indent For our next theorem, we wish to show that the dimension is continuous, given the right assumptions. Here, we will assume that the contraction maps are all similarities, and that the open set condition holds. Under these assumptions, the dimension depends on the similarity ratios of the iterated function system, so we first need to show that these ratios are continuous in the parameter.
\begin{prop}
Let $P$ and $X$ be metric spaces. Let $f: P \times X \to X$ be a
continuous family of similarities, with similarity ratio $r_p \in
(0,1)$. Then the function $r:P \to (0,1)$ is continuous. \end{prop}
\noindent \emph{Proof.} Let $\epsilon > 0$, and let $x,y \in X$,
with $x \ne y$. Then,
\[
d(f(p,x),f(p,y)) \leqslant d(f(p,x),f(q,x)) + d(f(q,x),f(q,y)) +
d(f(q,y),f(p,y)).
\] Similarly,
\[
d(f(q,x),f(q,y)) \leqslant d(f(q,x),f(p,x)) + d(f(p,x),f(p,y)) +
d(f(p,y),f(q,y). \] By subtracting, and combining the above two equations, we get
\[|d(f(p,x),f(p,y)) - d(f(q,x),f(q,y))| \leqslant d(f(p,x),f(q,x)) + d(f(q,y),f(p,y)). \]
Choose a $\delta > 0$ such that if $d(p,q) < \delta$, then
\[\max(d(f(p,x),f(q,x)),d(f(q,y),f(p,y))) < \frac{\epsilon}{2}d(x,y).\]
Then, if $d(p,q) < \delta$,
\[|d(f(p,x),f(p,y)) - d(f(q,x),f(q,y))| \leqslant \epsilon d(x,y).\]
But $f$ is a similarity, so
\[|r_pd(x,y) - r_qd(x,y)| < \epsilon d(x,y),\]
and therefore $|r_p - r_q| < \epsilon$. Thus $r:P \to (0,1)$ is continuous. $\hfill \square$
\\ \\
\indent Now that we know the similarity ratio is continuous in the parameter, we can easily show that the dimension is continuous. Since the dimension is implicitly defined by the ratios (as in theorem 1.12), the proof is primarily an application of the implicit function theorem.
\begin{thm}
Let $P$ and $X$ be metric spaces. Let $f_i: P \times X \to X$ be a
continuous iterated function system. Furthermore, let each of the
functions $f_i(p,\cdot)$ be similarities that satisfy the open set
condition for any $p \in P$. Let $A(p)$ be the invariant set for the
iterated function system for each $p$. Then the dimension of $A(p)$
is a continuous function of $p$. \end{thm} \noindent \emph{Proof.}
Let $g:(0,1)^n \times \mathbb{R} \to \mathbb{R}$ be given as $g(r,s)
= \sum_{i=1}^{n}r_i^s.$ Let $r \in (0,1)^n.$ Since $g(r,0)=n$ and
$\lim_{s \to \infty}g(r,s)=0,$ there must be some $s \in
\mathbb{R}_+$ such that $g(r,s) = 1$, by the intermediate value
theorem. Since $g$ is strictly decreasing in $s$, there is a unique
$s \in \mathbb{R}$ that satisfies this equation. Let the function
thus defined be called $D:(0,1)^n \to \mathbb{R}$. Since
\[\frac{\partial}{\partial s}\:g(r,s) = \sum_{i=1}^{n}(\log r_i)
r_i^s \ne 0\] for any s, the implicit function theorem implies this
function $D$ is continuous at $r$. But $r$ was arbitrary, so $D$ is
globally continuous. Now let $f_i$ be as above, and let $r_i(p)$ be
the similarity ratio for $f_i$ at $p$. Let $r(p) = (r_1(p),
r_2(p),...,r_n(p)).$ By proposition 2.4 , the function $r:P \to
(0,1)^n$ is continuous. Thus, the composition $D \circ r: P \to
\mathbb{R}$ is continuous. But since the iterated function systems
all consist of similarities that satisfy the open set condition,
theorem 1.12 says $D \circ r$ is exactly the dimension of the
invariant set. $\hfill \square$
\\ \\
\indent We now construct an example that shows that the Hausodorff
dimension of the invariant set need not be continuous if the open
set condition does not hold. That is, the open set condition is not
a superfluous assumption, and dimension can be discontinuous
otherwise.
\begin{ex}Let $P = [0,1],\: X = \mathbb{R},\: f_1(p,x) = \frac{x}{3} - p,\: f_2(p,x) = \frac{x}{3} + p.$ Clearly these are similarities from $\mathbb{R}$ into itself. For any $p \ne 0$, the open set condition holds, for $(-3p,3p)$ gets mapped to $(-2p,0)$ and $(0,2p)$. Thus the dimension of the invariant set is the solution to $2(\frac{1}{3})^s = 1$, which is $\frac{\log{2}}{\log{3}}$. However, when $p=0$, $f_1(0,x) = f_2(0,x)$. Then there is only a single point, $\{0\}$, which is invariant, and the dimension of one point is zero. Thus, the dimension as a function of $p$ is discontinuous at $0$.
\end{ex}

\end{document}